\documentstyle[12pt,amssymb]{article}

\newenvironment{diagram}{%
   \begingroup
   \arraycolsep0.1em
   
}{%
   \endgroup%
}
\newenvironment{commdiag}{%
   \begin{diagram}
   $$
   \begin{array}{cccccccccccccc}
}{
   \end{array}
   $$%
   \end{diagram}%
}
\newlength{\rightarrowlength} \rightarrowlength2em
\newlength{\leftarrowlength} \leftarrowlength2em
\newlength{\downarrowlength}  \downarrowlength4ex
\newlength{\diagarrowlength}  \diagarrowlength2em

\def\Rarrow{%
   \mathop{\makebox[\rightarrowlength]{\rightarrowfill}}\limits%
}
   
\def\Requal{%
   \mathop{%
   \rlap{\raise0.75ex\hbox{\makebox[\rightarrowlength]{\hrulefill}}}%
   \raise0.25ex\hbox{\makebox[\rightarrowlength]{\hrulefill}}}\limits%
}
\def\Darrow{%
   \left\downarrow\parbox{0cm}{\rule{0cm}{\downarrowlength}}\right.%
}

\def\Dhook{%
   \left\downarrow\parbox{0cm}{\rule{0cm}{\downarrowlength}}\right.%
   \kern-1.46ex\raise0.50\downarrowlength\hbox{$\scriptscriptstyle\cap$}%
}

\def\rlabel#1{%
   \rlap{$\scriptstyle#1$}
}
\def\llabel#1{%
   \llap{$\scriptstyle#1$}
}

\def\to{\longrightarrow}

\def\phi{\varphi}
\def\epsilon{\varepsilon}
\def\tilde{\widetilde}
\def\hat{\widehat}
\def\bar{\overline}
\def\({\left(}
\def\){\right)}
\def\O{{\cal O}}
\def\Z{{\cal Z}}
\def\I{{\cal I}}
\def\J{{\cal J}}
\def\W{{\cal W}}

\newtheorem{satz}{Satz}[section]
\newtheorem{theorem}[satz]{Theorem}

\newtheorem{remark}[satz]{Remark}
\newtheorem{lemma}[satz]{Lemma}

\newtheorem{example}[satz]{Example}
\newenvironment{varthm*}[1]{\trivlist\item[]\it{\bf #1.}}{\endtrivlist}
\def\parag#1{\par\addvspace{\bigskipamount}\noindent{\bf#1.}}

\def\startproof{\addvspace{\bigskipamount}\noindent}
\def\proof{\startproof{\it Proof. }}
\def\proofof#1{\startproof{\it Proof of #1.}}
\def\qedsymbol{\frame{\rule[0pt]{0pt}{8pt}\rule[0pt]{8pt}{0pt}}}
\def\qed{\nopagebreak\hspace*{\fill}\qedsymbol\par\addvspace{\bigskipamount}}

\def\bbP{{\Bbb P}}

\def\bbC{{\Bbb C}}
\def\bbF{{\Bbb F}}
\def\bbP{{\rm I\!P}}
\def\bbF{{\rm I\!F}}
\def\A{{\frak{a}}}
\def\be{\begin{eqnarray*}}
\def\ee{\end{eqnarray*}}
\def\lreqn#1#2{
   \begin{trivlist}\item[]
      $\begin{array}[t]{r} #1 \end{array}$ \\
      \hspace*{\fill}$\begin{array}[t]{@{}l} #2 \end{array}$
   \end{trivlist}
}

\def\tensor{\otimes}

\def\Bigwith{\ \Big\vert\ }
\def\liste#1#2#3{\mbox{$#1_{#2},\dots,#1_{#3}$}}
\def\eqdef{=_{\mrom{def}}}

\def\with{\mid}

\def\inverse{^{-1}}
\def\isom{\cong}

\def\rounddown#1{\left\lfloor#1\right\rfloor}
\def\mrom#1{{\rm #1}}

\def\Pic{\mathop{\mrom{Pic}}\nolimits}
\def\max{\mathop{\mrom{max}}\nolimits}
\def\ev{\mathop{\mrom{ev}}\nolimits}
\def\supp{\mathop{\mrom{supp}}\nolimits}

\def\m{{\frak m}}

\def\lra{\longrightarrow}
\def\calo{{\cal{O}}}
\def\mod{\mathop{\mrom{mod}}\nolimits}
\def\len{\mathop{\mrom{length}}\nolimits}
\def\length{\len}

\def\hat#1{\widehat{#1}}
\def\a{\alpha}
\def\e{\epsilon}
\def\b{\beta}

\begin{document}
\title{Cyclic coverings and higher order embeddings\\ of algebraic varieties}
\author{Th.\ Bauer,
  S.\ Di Rocco, 
  T.\ Szemberg\thanks{Supported by a research grant of Polish Academy of Sciences}}
\maketitle

%****************************************************************************

\begin{abstract}
%\noindent{1991 {\em Mathematics Subject Classification}. Primary ;
%   Secondary }\\
%\noindent{{\em Keywords and phrases.}}
  Let $Y$ be a cyclic covering of an algebraic variety $X$. Given
  a line bundle on $X$ we give criteria for its pull-back to $Y$
  to define a higher order embedding.
\end{abstract} 

%****************************************************************************

%\setcounter{section}{-1}
\section*{Introduction}

   In recent years there has been considerable interest in understanding
   under which circumstances linear series on algebraic varieties restrict
   surjectively to zero-dimensional subschemes and to collections of fat
   points.  Equivalently, one asks for the {\em order} of a given projective
   embedding of an algebraic variety.  The concepts of higher order
   embeddings are captured by the
   notions of $k$-very ampleness, introduced in \cite{BelFraSom89},
   and $k$-jet ampleness, studied in \cite{BelSom93b}
   and \cite{Dem}.
   More geometrically, an algebraic variety embedded in $\bbP^N$ via
   a $k-$very ample line bundle has no $(k+1)-$secant $(k-1)-$plane
   $\bbP^{k-1}\subset\bbP^N$. The embedding given by a $k-$jet
%   ample line bundle has stronger geometrical constrains related
   ample line bundle has stronger geometrical constraints related
   to the higher osculating planes.

   By now, the situation on surfaces is quite well understood, mainly
   thanks to the availability of powerful methods such as a Reider type
   criterion for $k$-very ampleness \cite{BelFraSom89} and the use of
   $\Bbb Q$-divisors in connection with the Kawamata-Viehweg vanishing
   theorem \cite{Laz93}.  In higher dimensions, however, the problem seems
   to be much more difficult.  While there are highly interesting general
   results on the separation of jets due to Demailly \cite{Dem}, there is
   still a lack of practical criteria that allow to determine the order
   of a given embedding.  The purpose of this paper is to study higher
   order embeddings of cyclic coverings $\pi:Y\lra X$, via line
   bundles given by pulling back "sufficiently positive" line
   bundles on $X$.
   Given a line bundle $L$ on $X$ we relate the
   order of the embedding defined by $\pi^*L$ to that of $L$ and certain
   rank $1$ summands of the vector bundle $L\tensor\pi_*\calo_Y$.
   The main results are expressed in the Theorems \ref{jet}
   and \ref{very}.

   In Section 1 we recall briefly the concepts of higher order
   embeddings. Section 2 is devoted to an exposition of the
   $k-$jet ample case. In Section 3 we prove our result concerning
   $k-$very ampleness. Although the rigorous proof is technically
   involved, we hope that the ideas behind are fairly transparent.
   Finally, the examples provided in the last section show that
   our results are sharp.

%\parag{\it Notation and Conventions}
   We work throughout over the field $\bbC$ of complex numbers.

\parag{\it Acknowledgments}
%   We would like to thank the Max-Planck-Institut f\"ur Mathematik, where
%   the first version of this paper was written, for financial support
%   which enabled us to come together and for providing excellent
%   working conditions.
   The first version of the paper was written during the authors' stay at
   the Max-Planck-Institut f\"ur Mathematik in Bonn.
   We would like to thank the Institute for financial support and for
   providing excellent working conditions. Our stay has been most fruitful
   and enjoyable.  
%****************************************************************************

\section{Higher order embeddings}

   We start by recalling the notions of $k$-jet ampleness and
   $k$-very ampleness
   that capture the concept of a {\em higher order embedding} in two
   different ways: the first notion requires the
   simultaneous separation of jets at finitely many points,
   whereas the second asks for the surjectivity of 
   the restriction to 0-dimensional
   subschemes of certain length.

   A line bundle $L$ on a smooth projective variety $X$
   is called {\em $k$-jet ample},
   if the evaluation map
   $$
      H^0(X,L) \to
      H^0\(X,L\tensor\O_X/\(\m_{y_1}^{k_1}\tensor\dots\tensor \m_{y_r}^{k_r}\)\)
   $$
   is surjective for any choice of distinct points $\liste y1r$ in $X$
   and positive integers $\liste k1r$ with $\sum k_i=k+1$.
   Further,
   $L$ is called {\em $k$-very ample}, if
   for every zero-dimensional subscheme $(\Z,\O_Z)\subset(X,\O_X)$
   of length $k+1$ the natural map
   $$
      H^0(X,L) \to H^0\(X,L\tensor\O_Z\)
   $$
   is surjective.

   As the definitions suggest, $k-$jet ampleness implies $k-$very
   ampleness \cite[2.2]{BelSom93b} and
   both notions are equivalent to global generation
   in the case $k=0$ and very ampleness in the case $k=1$.

In the case when $X=\bbP^2$ the notion of $k$-the order embedding is another
way of dealing with non ``superabundant embeddings". We refer to \cite{Ge} and \cite{BelSom96} for more details.
   
   For the purposes of the proof of Theorem \ref{very} below, it will
   be useful to have a slightly more general definition of $k$-very
   ampleness:
   
A subspace $V\subset H^0(X,L)$ will be called $k$-very ample
   on an open subset $U\subset X$, if for every
   zero-dimensional subscheme $(\Z,\O_\Z)\subset(U,\O_U)$
   of length $k+1$ the natural map
   $$
      \ev_\Z: V \to H^0\(X,L\tensor\O_\Z\)
   $$
   is surjective.  It is then obvious that if $V\subset H^0(X,L)$ is
   $k$-very ample on $U$, and if $(\Z,\O_Z)\subset(U,\O_U)$
   is a 0-dimensional
   subscheme of length $\ell$, then the subspace
   $$
      \ev_\Z\inverse(g)\subset V
   $$
   is $(k-\ell)$-very ample on $U\setminus\supp(\Z)$ for every
   $g\in H^0(X,L\tensor\O_\Z)$.

%****************************************************************************

\section{Coverings and $k$-jet ampleness}

\begin{theorem}\label{jet}
   Let $X$ be a smooth projective
   variety and $B\subset X$ a smooth divisor.
   Let $M$ be a line bundle on $X$ such that
   $\O_X(dM)\isom\calo_X(B)$ and let $\pi:Y\lra X$
   be the cyclic covering of degree $d$ defined by $M$.
   Let $L$ be a line bundle on $X$ and $k$ a non-negative
   integer.
   If $L-qM$ is $(k-q)$-jet ample for $q=0,\ldots,\min(k,d-1)$,
   then $\pi^*L$ is $k$-jet ample.
\end{theorem}

\proof
   Let $\bar{M}$ be the total space of $M$ and
   $p:\bar{M}\lra X$ the bundle projection.
   Let $s_B\in H^0(X,\calo_X(B))$ be a section
   whose divisor of zeros is $B$, and let $\tau\in H^0(\bar{M},p^*M)$
   be the tautological section. As usual, $Y$ may be viewed
   as the divisor of zeros of the section $p^*s_B-\tau^d$.
   Letting $t$ be the restriction of $\tau$ to $Y$
   and $\pi:Y\to X$ the restriction of $p$,
   the projection formula gives the decomposition
   \begin{equation}\label{decomposition}
     H^0(Y,\pi^*L)=\bigoplus_{q=0}^{d-1}t^q\pi^*H^0(X,L-qM)
   \end{equation}
   corresponding to the eigen-values of the action of the
   primitive covering automorphism $\varphi$.

   Let $y_1,\ldots,y_r\in Y$ be points and $k_1,\ldots,k_r$
   positive integers with $\sum k_i=k+1$. Given a simultaneous jet
   $$
      J\in
      \bigoplus_{i=1}^rH^0(Y,\pi^*L\tensor\calo_Y/\m_{y_i}^{k_i})=
      H^0(Y,\pi^*L\tensor\calo_Y/\m_{y_1}^{k_1}\tensor\ldots
      \tensor\m_{y_r}^{k_r}) \ ,
   $$
   we decompose $J$ into a sum of simultaneous jets $J=J_1+\ldots+J_r$,
   where $J_i=(j_{1i},\ldots,j_{ii},\ldots,j_{ri})$, with
   $j_{li}$ the zero jet of order $k_l$ at $y_l$ for $l\neq i$.
   It is enough to find  sections $s_i\in H^0(Y,\pi^*L)$ such that
   $$
      s_i \mod \m_{y_1}^{k_1}\tensor\ldots\tensor\m_{y_r}^{k_r}=J_i
      \quad\mbox{ for } i=1,\ldots,r \ ,
   $$
   since then $s=s_1+\ldots+s_r$ has the desired jet $J$.
   We may therefore assume $J=J_1$.
   Now, the idea is to construct a section $s$ explicitly out
   of sections in the line bundles $L-qM$, $q\geq 0$. We distinguish
   three cases.

   {\em Case 1.} Suppose that $y_1\notin R$ and that none of the
   points $y_2,\ldots,y_r$ lies in the orbit $\pi^{-1}(\pi(y_1))$.
   This is the easiest case. The desired sections is obtained
   as a pull-back of a section in $L$. More precisely
   we can view the jet $j_{11}$ as a jet at $\pi(y-1)$,
   via the isomorphism of local
   rings $\calo_{Y,y_1}\isom\calo_{X,\pi(y_1)}$ induced by $\pi$.
   Since $L$ is $k-$jet ample, there
   is a section $s\in H^0(X,L)$ such that
   $$
      s \mod \m_{\pi(y_1)}^{k_1} = j_{11}
   $$
   and
   $$
      s \mod \m_{\pi(y_i)}^{k_i} = 0 \quad\mbox{ for } i=2,\ldots,r
      \ .
   $$
   Then of course $\pi^*s \mod
   \m_{y_1}^{k_1}\tensor\ldots\tensor\m_{y_r}^{k_r} = J_1$.

   {\em Case 2.} Suppose that $y_1,\ldots,y_l$ (with $l\ge 2$) lie in the orbit
   $\pi\inverse(\pi(y_1))$ and that none of the points $\liste y{l+1}r$ does.
   This case is more difficult since we have to separate points in a fiber
   of the covering. The construction of $s$ builds upon the cyclic group
   at hand.
   Let $k'=\max(\liste k1l)$. Observe that 
   $k'\le k+1/(l-1)-\sum_{i=l+1}^r k_i$.
   Let $\hat j_{11}$ be any preimage of the jet $j_{11}$ under the surjective
   map 
   $$
      \O/\m_{y_1}^{k'}\to\O/\m_{y_1}^{k_1} \ .
   $$
   Using now again the identification $\O_{Y,y_1}\isom\O_{X,\pi(y_1)}$ as
   well as the facts that $t(y_1)\ne 0$ and that the line bundle $L-qM$
   is $(k-q)$-jet ample for $0\le q\le l-1$, we can find sections $\liste s1l$
   such that
   \be
      \pi^*s_1 \mod \m_{y_1}^{k'} &=& \alpha_1\cdot \hat j_{11} \\   
      t\cdot\pi^*s_2 \mod \m_{y_1}^{k'} &=& \alpha_2\cdot \hat j_{11} \\   
      \vdots \\
      t^{l-1}\cdot\pi^*s_l \mod \m_{y_1}^{k'} &=& \alpha_l\cdot \hat j_{11} 
   \ee
   and
   $$
      t^q\cdot\pi^*s_{q+1} \mod \m_{y_i}^{k_i} = 0 
      \quad\mbox{ for } 0\le q\le l-1 \mbox{ and } l+1\le i \le r \ .   
   $$
   for any given complex numbers $\liste{\alpha}1l$.
   Now, let $(\liste{\alpha}1l)$ be a solution of the system of linear
   equations
   \begin{equation}\label{linsys1}
   \begin{array}{rcl}
   \alpha_1+\ldots+\alpha_l&=&1\\
   \alpha_1+\epsilon^{\beta_2}\alpha_2+\ldots+\epsilon^{(l-1)\beta_2}\alpha_l&=&0\\
   &\vdots& \\
   \alpha_1+\epsilon^{\beta_l}\alpha_2+\ldots+\epsilon^{(l-1)\beta_l}\alpha_l&=&0
   \end{array}
   \end{equation}
   where $\beta_2,\ldots,\beta_l$ are integers such that
   $y_i=\varphi^{\beta_i}(y_1)$.
   (Note that the determinant of the system (\ref{linsys1})
   is just the Vandermonde determinant of the numbers
   $1,\epsilon^{\beta_2},\ldots,\epsilon^{\beta_l}$.)
   Consider the section
   $s=\pi^*s_1+\ldots+t^{l-1}\pi^*s_l$.
   We have $s \mod \m^{k'}_{y_1} = \hat j_{11}$, which implies
   $$
      s \mod \m_{y_1}^{k_1}=j_{11} \ ,
   $$
   as required. Moreover we have
   \be
      s\mod\m_{y_i}^{k'}&=&(\varphi^{\beta_i})^*s\mod\m_{y_1}^{k'}\\
      &=&
      \pi^*s_1+\epsilon^{\beta_i}t\pi^*s_2+\ldots+\epsilon^{(l-1)\beta_i}
      t^{l-1}\pi^*s_l \mod\m_{y_1}^{k'}\\
      &=&
      (\alpha_1+\epsilon^{\beta_i}\alpha_2+\ldots+
      {(\epsilon^{\beta_i})}^{l-1}\alpha_l)\hat j_{11}=0 \ ,
   \ee
   so that $s\mod\m_{y_i}^{k_i}=0$ for $i=2,\ldots,l$,
   and of course
   $$
      s\mod\m_{y_i}^{k_i}=0 \mbox{ for } i=l+1,\ldots,r \ .
   $$

   {\em Case 3.}
   Suppose $y_1\in R$. Since $B$ (and hence $R$) is smooth, there
   are local coordinates $u_1,\ldots,u_n$ at $y_1$ and
   $v_1,\ldots,v_n$ at $\pi(y_1)$ such that $\pi$ is locally given as
   $$
      \pi(u_1,\ldots,u_n)=(u_1^d,u_2,\ldots,u_n)=(v_1,\ldots,v_n)
      \ .
   $$
   In these coordinates $t$ is given by $u_1$ and $j_{11}$ may be
   written as
   $$
      j_{11}=\sum_{i_1+\ldots+i_n\leq k_1}a_{i_1,\ldots,i_n}
   u_1^{i_1}\cdot\ldots\cdot u_n^{i_n}=
   \sum_{l=0}^{d-1}u_1^{l}
   \sum_{\begin{array}{c}\scriptstyle i_1\equiv l\mod d\\
   \scriptstyle i_1+\ldots+i_n\leq k_1\end{array}}
   a_{i_1\ldots i_n}{(u_1^d)}^{\frac{i_1-l}{d}}u_2^{i_2}\cdot\ldots
   \cdot u_n^{i_n}.
   $$
   This splitting reflects the fact that only jets containing
   powers of $u_1$ divisible by $d$ arise as pull-backs of jets on $X$.
   Since $L-qM$ is $(k-q)$-jet ample for $q=0,\ldots,\min(d-1,k)$,
   there are sections $s_q\in H^0(X,L-qM)$ such that
   $$s_q\mod\m_{\pi(y_1)}^{k_1-q}=
   \sum_{\begin{array}{c}\scriptstyle i_1\equiv l\mod d\\
         \scriptstyle\frac{i_1-q}{d}+i_2+\ldots+i_n\leq k_1-q\end{array}}
   a_{i_1\ldots i_n}v_1^{\frac{i_1-l}{d}}v_2^{i_2}\cdot\ldots
   \cdot v_n^{i_n}\ .
   $$
   and
   $$
      s_q\mod\m_{\pi(y_i)}^{k_i}=0 \quad\mbox{ for } i=2,\ldots,r\ .
   $$
   It is now easy to check that the section
   $$
      s=\sum_{q\geq 0}t^q\pi^*s_q
   $$
   has the prescribed jets.
\qed

   Working as in cases 1 and 2 of the above proof, one gets immediately
   the same result for unbranched coverings:

\begin{theorem}(Unbranched case)
   Let $\pi:Y\lra X$ be an unbranched cyclic covering of degree $d$,
   defined by a line bundle $M$ with $\O_X(dM)=\O_X$.
   Let $L$ be a line bundle on $X$ and $k$ a non-negative integer.
   If $L-qM$ is $(k-q)$-jet ample for $q=0,\ldots,\min(k,d-1)$,
   then $\pi^*L$ is $k$-jet ample.
\end{theorem}

%****************************************************************************

\section{Coverings and $k$-very ampleness}

%   In this section ...

   In Theorem \ref{very} below we will assert the $k$-very ampleness
   of a pullback of an ample line bundle $L$ under suitable hypotheses
   on the line bundles $L-qM$, $q\ge 0$.
   In order to get a useful criterion, one needs to formulate
   a delicate numerical hypothesis.
   To this end, given positive integers $k$ and $\ell$,
   we introduce the
   abbreviations
   $\gamma(k,\ell)=1$ if $k/ \ell$ is an integer, and
   $\gamma(k,\ell)=0$ otherwise, and we set
   $$
      \tau(k,\ell) = k-\rounddown{\frac k\ell}-\ell+\gamma(k,\ell)+1 \ .
   $$
   Further, we use
   the notation
   $$
      \sigma(k,d,q)=
            \max\left\{\tau(k+1,\ell)
               \Bigwith q+1\le\ell\le\min(d,k+1)\right\}-1
   $$
   if $q>0$ and $\sigma(k,d,0)=k$.

   Our result can then be stated as follows:

\begin{theorem}\label{very}
   Let $X$ be a smooth projective
   variety and $B\subset X$ a smooth divisor.
   Let $M$ be a line bundle on $X$ such that
   $\O_X(dM)\isom\calo_X(B)$ and let $\pi:Y\lra X$
   be the cyclic covering of degree $d$ defined by $M$.
   Let $L$ be a line bundle on $X$ and $k$ a non-negative
   integer.
   If $L-qM$ is $\sigma(k,d,q)$-very ample for $q=0,\ldots,\min(k,d-1)$
%%%   and if $L$ is $k$-very ample
   then $\pi^*L$ is $k$-very ample.
\end{theorem}

   One can prove the $k$-very ampleness of $\pi^*L$ under
   the assumption that $L-qM$ is $(k-q)$-very ample for
   $0\le q\le\min(k,d+1)$ by arguments very similar to the ones used in
   the proof of Theorem \ref{jet}. The point however here is that
   $k$-very ampleness holds already under the much weaker (but also more
   intricate) hypotheses involving the numbers $\sigma(k,d,q)$.
   For instance, if $L$ is 2-very ample, then $\pi^*L$ is also 2-very
   ample as soon as both $L-M$ and $L-2M$ are globally generated.
   If $\pi$ is of degree $2$ it is enough to check that $L-M$
   is globally generated.
   We provide a list of explicit values of
   $\sigma(k,d,q)$ in Remark \ref{sigma values}.

   In order to simplify the exposition of the proof of the theorem,
   we start by stating two lemmas.
   We begin with an elementary observation on local
   rings:

\begin{lemma}\label{alg}
   Let $(\O,\m)$ be a local ring, and let $\liste I1\ell$, $\ell\ge 2$,
   be ideals contained in $\m$ such that
   $\length\O/I_1\ge\length\O/I_i$ for all $i$.
   Let $k=\sum_{i=1}^{\ell}\len\O/ I_i$.
   If $k<\infty$, then
   $$
      \length\O/ I_2\cap\dots\cap I_\ell\le
      \tau(k,\ell) \ .
   $$
\end{lemma}

\proof
   We have the inequality
   \be
      \length\O/ I_2\cap\dots\cap I_\ell
      &\le&\sum_{i=1}^{\ell}\len\O/ I_i
      -\len\O / I_2+(I_3\cap\dots\cap I_\ell) \\
      &&-\dots
      -\len\O / I_2+\dots+I_\ell
      \le k-\len\O/ I_1-(\ell-2) \ .
   \ee
   The assertion follows now from  the fact that
   $\len\O/ I_1\ge\rounddown{k/ \ell}+1-\gamma(k,\ell)$.
\qed

   We will also need the following

\begin{lemma}\label{num}
   Let $\liste\ell1m$, $\liste K1r$ and $q$ be positive integers,
   and let $K=\sum_{i=1}^r K_i$.
   If $\ell_i\ge 2$ for all $i$, then
   $$
      \sum_{i=1}^m \tau(K_i,\ell_i)
      +\sum_{i=m+1}^r \rounddown{\frac{K_i}{q+1}}
      \le
      \max\left\{\tau(K,\ell)
      \Bigwith
         \min(\liste\ell1q,q+1)\le\ell\right\} \ .
   $$
\iffalse
   \lreqn{
      \displaystyle
      \sum_{i=1}^m
         \(K_i-\rounddown{\frac{K_i}{\ell_i}}-\ell_i+1+\gamma(\ell_i,K_i)\)
      +\sum_{i=m+1}^r \rounddown{\frac{K_i}{\ell_i}}
    }{\le
      \displaystyle
      \max\left\{K-\rounddown{\frac K\ell}-\ell+1+\gamma(\ell,K)\Bigwith
         \min(\liste\ell1q,q+1)\le\ell<\infty\right\} \ .
   }
\fi
\end{lemma}

\proof
  First we note that it is enough to prove
  \begin{equation}\label{first}
      \sum_{i=1}^m \tau(K_i,\ell_i)
      +\sum_{i=m+1}^r \rounddown{\frac{K_i}{q+1}}
      \le
      \tau(K,\ell)
  \end{equation}
  with $\ell=\max(\ell_1,\dots,\ell_m)$.
  Due to the fact that
   $$
      \rounddown{\frac{K_i}{q+1}}\le
      K_i-\rounddown{\frac{K_i}{2}}-1+\gamma(K_i,2)=\tau(K_i,2) \ .
   $$
   we may assume $r=m$ (by setting $l_i=2$ for $m<i\le r$).
   Then (\ref{first}) reads
   $\sum_{i=1}^r \tau(K_i,\ell_i)\leq \tau(K,\ell)$
   and it suffices to prove it for $r=2$. We may also
   assume $\ell_1\leq \ell_2=\ell$.
   From
   $$\rounddown{\frac{K_1}{\ell_1}}\geq\rounddown{\frac{K_1}\ell}+
     \gamma(K_1,\ell_1)\cdot\min(1,\ell-\ell_1)$$
   and
   $$\rounddown{\frac{K_1}{\ell}}+\rounddown{\frac{K_2}\ell}\geq
     \rounddown{\frac{K}\ell}-1+\gamma(K_1,\ell)+\gamma(K_2,\ell)-
     \gamma(K_1,\ell)\cdot\gamma(K_2,\ell)$$
   we get
   \begin{center}
   \begin{tabular}{rcl}
   $\tau(K_1,\ell_1)+\tau(K_2,\ell_2)$&$\leq$&
     $K-\rounddown{\frac{K}\ell}-\ell+\gamma(K_1,\ell)\cdot\gamma(K_2,\ell)+1+$\\
     &&$+[2-\ell_1]+[\gamma(K_1,\ell_1)-\gamma(K,\ell)-\gamma(K_1,\ell_1)\cdot
     \min(1,\ell-\ell_1)].$
   \end{tabular}
   \end{center}
   Now, the assertion follows by observing that
   $\gamma(K_1,\ell)\cdot\gamma(K_2,\ell)\leq\gamma(K,\ell)$
   and the terms in the square brackets are non-positive.
\qed

\proofof{Theorem \ref{very}}
   Let $({\Z},{\O}_{{\cal Z}})\subset Y$ be a zero-dimensional subscheme of
   length $k+1$, defined by the ideal $\I_Z$, and
   supported on the points $y_1,\dots,y_r$.  Let $k_i=\len(\O_{\Z,y_i})$.
   Given an element
   $$g=(g_1,\dots,g_r)\in\bigoplus_{i=1}^{r} H^0(\pi^*L\tensor\O_{\Z,y_i})
   =H^0(\pi^*L\tensor\O_{\Z}) \ ,$$
   we will construct a section $s\in H^0(Y,\pi^*L)$ whose image
   in $H^0(\pi^*L\tensor\O_{\Z})$ is $g$.

   The strategy is different to that in the proof of Theorem \ref{jet}.
   Roughly speaking, for $q\geq 0$, we will approximate $s$ by building
   a sequence of affine subspaces
   $$\dots \subset V^q_2\subset V^q_1\subset V^q_0\subset
     H^0(X,L-qM)$$
   consisting of sections $s_{q,i}$ such that the length
   of the subscheme of $\Z$ where $\sum t^q\pi^*s_{q,i} \mod \I_\Z$
   agrees with $g$ increases with $i$
   reaching $\length\Z$ in the last step.
   First of all we take care of subschemes of $\Z$ supported
   in the ramification locus of $\pi$.
   Starting with a fiber containing most of the points $y_i$
   we deal then with subschemes of $\Z$ supported in regular fibers of $\pi$.
   Each step consists now of fixing sections in $V^q_i$ for
%   $q$ appropriatly big and imposing new conditions on sections
   $q$ appropriately big and imposing new conditions on sections
   in $V^q_i$ for lower $q$ thus defining subspaces $V^q_{i+1}$.
   Positivity assumptions of the theorem assure that we never
   run out of sections. Let us now turn to the details.

   After reordering $y_1,\dots,y_r$ we may assume that
   \be
   y_1,\dots,y_{n_1-1}&\in& R \\
   y_{n_1},\dots,y_{n_2-1}&\in&\pi^{-1}(x_1)\\
   &\vdots \\
   y_{n_m},\dots,y_{n_{m+1}-1}=y_r&\in&\pi^{-1}(x_m) \ .
   \ee
   for distinct points $x_1,\dots,x_n\in X\setminus B$.
   Letting $l_i=n_{i+1}-n_i$, we may also assume that
   $l_1\leq\dots\leq l_m$ and $k_{n_i}=\max(k_{n_i},\dots, k_{n_{i+1}-1})$.
   We will use the abbreviations
   $\I_i=\I_{\Z,y_i}$ for $i=1,\dots, r$ and
   $K_i=k_{n_1}+\dots+k_{n_{i+1}-1}$ for $i=1,\dots, m$.

   Identifying the local rings of the points in the same fiber by means of $\pi$,
   we consider the ideals
   $$\J_i\eqdef\I_{n_i+1}\cap\dots\cap\I_{n_{i+1}-1}\subset\O_{Y,y_{n_i}}.$$
   and for $i=1,\dots,m$ we denote by
   $$\rho_i:\O_{Y,y_{n_i}}/(\I_{n_i}\cap\J_i)\to\O_{Y,y_{n_i}}/ \J_i$$
   and
   $$ \eta_i:\O_{Y,y_{n_i}}/(\I_{n_i}\cap\J_i)\to\O_{Y,y_{n_i}}/\I_{n_i}$$
   the quotient maps. Applying Lemma \ref{alg} we get the inequalities
   \begin{equation}\label{R}
      \len({\O_{X,x_i}}/ \J_i)\leq\tau(K_i,l_i)
%%%   K_i-\rounddown{\frac{K_i}{l_i}}+\gamma(K_i,l_i)-l_i+1
   \end{equation}
   for all $i\in\{1,\dots,m\}$ such that $l_i\geq 2$.

   Let $y_i$ be one of the points $y_1,\dots,y_{n_1-1}$ on the ramification
   divisor $R$ and let $\W$ be the restriction of $\Z$ to $y_i$.
%   \O_{\W_i}=\O_{\Z,y_i}$ and  $\I_{\W_i}=\I_{\Z,y_i}$.
   Choosing coordinates $u_1,\dots, u_n$ at the point $y_i$ as
   in {\em Case 3} of the proof of Theorem \ref{jet}, we have a decomposition
   $$\I_{\W_i}=\bigoplus_{q=0}^{d-1}u^q\pi^{-1}\A_{i,q}\ ,$$
   where $\A_{i,q}$ are ideals in $\O_{X,\pi(y_i)}$,
   hence we get an isomorphism
   \begin{equation}\label{iso}
   \O_{\W_i}\cong \bigoplus_{q=0}^{d-1} \O_X/\A_{i,q}
   \end{equation}
   Note that $\A_{i,q}\subset\A_{i,q+1}$ for $q=0,\dots,d-2$, giving a cofiltration
   $$ \O_X/\A_{i,0}\twoheadrightarrow\O_X/\A_{i,1}\twoheadrightarrow\dots\twoheadrightarrow
   \O_X/\A_{i,d-1}\rightarrow 0 \ . $$
   This implies $\len(\O_X/\A_{i,q})\leq \rounddown{\frac{\len(\O_{\W_i})}{q+1}}$.
   Under the isomorphism (\ref{iso}) the element
   $g_i$ corresponds to a $d$-tuple $(g_{i,0},\dots,g_{i,d-1})$ for $i=1,\dots, n_1-1$.
   Since $L-qM$ is $\sigma(k,d,q)$-very ample, the subspace
   $$V_0^q=\{s\in H^0(X,L-qM) \with s\mod \A_{i,q}=g_{i,q}\mbox{ for }i=1,\dots,n_1-1\}$$
   is at least $\sigma_0(k,d,q)$-very ample on $X\setminus R$, where
   $$\sigma_0(k,d,q)=\sigma(k,d,q)-\rounddown{\frac{\len(\O_{\W_i})}{q+1}}$$
   for $q=0,\dots,\min(k,d-1)$. For any section $s$ of the form
   $$s=\pi^*s_0+t\pi^*s_1+\dots+t^{\min(k,d-1)}\pi^*s_{\min(k,d-1)}$$
   with $s_q\in V_0^q$ we then have
   $$s \mod \I_i = g_i \;\;\;\;\mbox{  in }\O_{Y,y_i}$$
   for $i=1,\dots,n_1-1$.

   From now on we proceed in $m$ steps.
   The recursive procedure works so that the $i$-th step takes care of
   what we generate in the $(m-i+1)$-st fiber.

   {\em Step 1.} For all $q\geq l_m$ such that $V_0^q \neq H^0(X,L-qM)$
   we fix a section $s_q\in H^0(X,L-qM)$. For $i=1,\dots,l_m$ we denote
   by $h_i$ the image $$(\phi^{\beta(n_m,n_m+i-1)})^*(\sum_{q\leq l_m} t^q
   \pi^*(s_q)) \mod \I_{n_m}\cap\dots\cap\I_{n_{m+1}-1}$$ in the local ring
   $\O_{Y,y_{n_m}}$.
   We choose preimages $\hat g_i$ of $g_{n_m+i-1}$ in $\O_{Y,y_{n_m}}/
   \I_{n_m}\cap\dots\cap \I_{n_{m+1}-1}$. For $q\leq l_m-1$ we let $V_1^q$
   be the subspace of $V_0^q$ of sections $s_q$ satisfying the following
   conditions in $\O_{Y,y_{n_m}}$:
   \be
      t^{l_m-1}\pi^* s_{l_m-1} \mod \J_m&=&\a_{l_m-1,1}\rho_m(\hat g_1-h_1)+\dots +\a_{l_m-1,l_m}\rho_m(\hat g_{l_m}-h_{l_m})\\
      &\vdots \\
      t\pi^* s_1 \mod \J_m&=&\a_{1,1}\rho_m(\hat g_1-h_1)+\dots +\a_{1,l_m}\rho_m(\hat g_{l_m}-h_{l_m})\\
      \pi^* s_0 \mod \I_{n_m}\cap \J_m&=&\a_{0,1}(\hat g_1-h_1)+\dots +\a_{0,l_m}(\hat g_{l_m}-h_{l_m}) \ ,
   \ee
   where $\a_{i,j}$ are solutions of the Vandermonde type systems of linear equations
   \be
   \a_{0,i}+\dots +\a_{l_m-1,i}&=&\delta_{1i} \\
   \a_{0,i}+\e^{\b(n_m,n_m+1)}\a_{1,i}+\dots +\e^{\b(n_m,n_m+1)}\a_{l_m-1,i}&=&\delta_{2i}\\
   &\vdots&\\
   \a_{0,i}+\e^{(l_m-1)\b(n_m,n_{m+1}-1)}\a_{1,i}+\dots +\e^{(l_m-1)\b(n_m,n_{m+1}-1)}\a_{l_m-1,i}&=&\delta_{l_m i} \ .
   \ee
   We observe that (\ref{R}) implies that $V_1^q$ remains at least
   $\sigma_1(k,d,q)$-very ample, where for $q\geq 1$
   $$
      \sigma_1(k,d,q)=
      \sigma_0(k,d,q)-\tau(K_m,l_m)
%%%   (K_m-\rounddown{\frac{K_m}{l_m}}-l_m+\gamma(l_m,K_m)+1)
   $$
   and for $q=0$
   $$\sigma_1(k,d,0)=\sigma_0(k,d,0)-\sum_{i=n_m}^{n_{m+1}-1} K_i=\sum_{i=1}^{n_m-1} k_i.$$
   Let $s$ be of the form $\pi^*s_0+\dots + t^l\pi^* s_l$, where $s_q\in V_1^q$, for $q\leq l_m-1$.
   Then in the local ring $\O_{Y,y_{n_m+i-1}}$ we have
   \be
   s \mod \I_{n_m+i-1}& = &(\phi^{\b(n_m,n_m+i-1)})^* s\mod \I_{n_m}   \\
   &=&\pi^*s_0+\e^{\b(n_m,n_m+i-1)}t\pi^*s_1+\dots+\e^{(l_m-1)\b(n_m,n_m+i-1)}t^{l_m-1}\pi^*s_{l_m-1}+\\
   &&+(\phi^{\b(n_m,n_m+i-1)})^* (\sum_{q\geq l_m} s_q)\mod \I_1=     \\
   &=& \a_{0,1}\eta_m (\hat g_1-h_1)+\dots
       +\a_{0,l_m}\eta_m(\hat g_{l_m}-h_{l_m})+\dots +   \\
   &&+\e^{\b(n_m,n_m+i-1)}
     \( \a_{1,1}\eta_m (\hat g_1-h_1)+\dots+
     \a_{1,l_m}\eta_m(\hat g_{l_m}-h_{l_m})\)+     \\
   &\vdots&                \\
   &&+\e^{(l_m-1)\b(n_m,n_m+i-1)}( \a_{l_m-1,1}\eta_m(\hat g_1-h_1)+
     \dots+ \a_{l_m-1,l_m}\eta_m(\hat g_{l_m}-h_{l_m})+   \\
   &&+\eta_m(h_i)=          \\
   &=&(\a_{0,1}+\e^{\b(n_m,n_{m}+i-1)}\a_{1,1}+\dots
      +\e^{(l_m-1)\b(n_m,n_{m}+i-1)}\a_{l_m-1,1})
      \eta_m(\hat g_1-h_1)+  \\
   &\vdots&     \\
   &&+(\a_{0,l_m}+\e^{\b(n_m,n_{m}+i-1)}\a_{1,l_m}+\dots
     +\e^{(l_m-1)\b(n_m,n_{m}+i-1)}\a_{l_m-1,l_m})
   \eta_m(\hat g_{l_m}-h_{l_m})+   \\
   &&+\eta_m(h_i)\\
   &=&g_{n_m+i-1}
   \ee

   {\em Step i.}
   Proceeding as in {\em Step 1}, we fix sections $s_q\in V_{i-1}^q$
   for $l_{m-i+1}\leq q < l_{m-i+2}$ and we construct subspaces
   $$V_i^q\subset V_{i-1}^q$$
   for $q<l_{m-i+1}$, such that for all $s_q\in V_i^q$ one has
   $$\sum_{q=0}^l t^q \pi^* s_q \mod \I_j = g_j$$
   for all $n_{m-i+1}\leq j< n_{m-i+2}$.
   Then $V_i^q$ is still $\sigma_i(k,d,q)$-very ample, where by (\ref{R})
   $$
   \sigma_i(k,d,q)=\sigma_{i-1}(k,d,q)-\tau(K_i,l_i)=
%%%   (K_i-\rounddown{\frac{K_i}{l_i}}-l_i+\gamma(l_i,K_i)+1)\\
   \sigma_{0}(k,d,q)-\sum_{\nu=i+1}^m \tau(K_\nu,l_\nu)
%%%   (K_{\nu}-\rounddown{\frac{K_{\nu}}{l_{\nu}}}-l_{\nu}+\gamma(l_{\nu},K_{\nu})+1)
   $$
   for $q\geq 1$ and
   $$\sigma_i(k,d,0)=\sigma_0(k,d,0)-K_i=
     \sum_{\nu=1}^{n_{m-i+1}-1} k_{\nu}.$$
   After $m-1$ steps, by Lemma \ref{num}, the space $V_{m-1}^q$ will still be $\sigma_{m-1}(k,d,q)$-very ample with
   $$\sigma_{m-1}(k,d,q)\geq\tau(K_1,l_1)
%%%   K_1-\rounddown{\frac{K_1}{l_1}}-l_1+\gamma(l_1,K_1)+1
   $$
   for $1\leq q<l_1$ and
   $$\sigma_{m-1}(k,d,0)\geq K_1.$$
   This guarantees that we can find sections $s_1,\dots,s_{l_1}$ such that
   $$ s \eqdef \sum_{q=0}^l t^q \pi^* s_q$$
   satisfies
   $$s\mod \I_{\Z} = g.$$
   This completes the proof of the theorem.
%%%   Step $m-1$ then completes the proof.   %%% m-1 ?
\qed

   As in section 2 we get immediately the same result for
   unbranched coverings. 
\begin{theorem}\label{veryunbranched}
   Let $X$ be a smooth projective
   variety.
   Let $M$ be a line bundle on $X$ such that
   $\O_X(dM)\isom\calo_X$ and let $\pi:Y\lra X$
   be the cyclic covering of degree $d$ defined by $M$.
   Let $L$ be a line bundle on $X$ and $k$ a non-negative
   integer.
   If $L-qM$ is $\sigma(k,d,q)$-very ample for $q=0,\ldots,\min(k,d-1)$
%%%   and if $L$ is $k$-very ample
   then $\pi^*L$ is $k$-very ample.
\end{theorem}
   
\begin{remark}\label{sigma values}\rm
   In order to convey some feeling for the numbers $\sigma(k,d,q)$
   appearing in the hypotheses of Theorem \ref{very}, we include a table
   for $d=15$ and $0\le k\le 15$. (The $k$-th column lists the values
   $\sigma(k,15,q)$ for $1\le q\le k$.)  For instance $\pi^*L$ is $2$-very
   ample if $L$ is so and both $L-M$ and $L-2M$ are globally generated,
   and it is $4$-very
   ample if $L$ is so and $L-M$ and $L-2M$
   are very ample and $L-3M$ and $L-4M$ are
   globally generated.

   \begin{center}\footnotesize
   \begin{tabular}{lrrrrrrrrrrrrrrrrrrrrrrrrrrrrrrrrrrrrrrrrrrrrrrrr}
             k&   0&   1&   2&   3&   4&   5&   6&   7&   8&   9&  10&  11&  12&  13&  14&  15& \\
      \hline
      $L- 1M$ &    &   0&   0&   1&   1&   2&   2&   3&   4&   4&   5&   6&   6&   7&   8&   9&\\
      $L- 2M$ &    &    &   0&   0&   1&   2&   2&   3&   4&   4&   5&   6&   6&   7&   8&   9&\\
      $L- 3M$ &    &    &    &   0&   0&   1&   2&   3&   3&   4&   5&   6&   6&   7&   8&   9&\\
      $L- 4M$ &    &    &    &    &   0&   0&   1&   2&   3&   4&   4&   5&   6&   7&   8&   8&\\
      $L- 5M$ &    &    &    &    &    &   0&   0&   1&   2&   3&   4&   5&   5&   6&   7&   8&\\
      $L- 6M$ &    &    &    &    &    &    &   0&   0&   1&   2&   3&   4&   5&   6&   6&   7&\\
      $L- 7M$ &    &    &    &    &    &    &    &   0&   0&   1&   2&   3&   4&   5&   6&   7&\\
      $L- 8M$ &    &    &    &    &    &    &    &    &   0&   0&   1&   2&   3&   4&   5&   6&\\
      $L- 9M$ &    &    &    &    &    &    &    &    &    &   0&   0&   1&   2&   3&   4&   5&\\
      $L-10M$ &    &    &    &    &    &    &    &    &    &    &   0&   0&   1&   2&   3&   4&\\
      $L-11M$ &    &    &    &    &    &    &    &    &    &    &    &   0&   0&   1&   2&   3&\\
      $L-12M$ &    &    &    &    &    &    &    &    &    &    &    &    &   0&   0&   1&   2&\\
      $L-13M$ &    &    &    &    &    &    &    &    &    &    &    &    &    &   0&   0&   1&\\
      $L-14M$ &    &    &    &    &    &    &    &    &    &    &    &    &    &    &   0&   0&\\
   \end{tabular}
   \end{center}
\end{remark}

%****************************************************************************

\section{Examples and applications}

   Theorems \ref{jet} and \ref{very} state, roughly speaking,
   that under suitable assumptions on the bundles $L-qM$ the
   positivity of a pullback $\pi^*L$ is at least as high as the
   positivity of $L$.
   It is then natural to ask whether the converse statement is also
   true. The following example shows that this is not the case:

\begin{example}\rm
   Let $(X,L)$ be a principally polarized abelian variety. A
   $d$-torsion point $M\in\Pic^0(X)$ give rise to a cyclic covering
   $\pi:Y\to X$ such that $(Y,\pi^*L)$ is an abelian variety of type
   $(1,\dots,1,d)$.  If $d$ is sufficiently large, then $\pi^*L$ will
   be very ample. (One knows by work of Debarre-Hulek-Spandaw
   \cite{DebHulSpa94} that it
   is enough to take $d>2^g$.)  But $L$ is not even globally
   generated.
   In the surface case, where the generation of jets
   and the $k-$very ampleness
   is well
   understood (see \cite{BauSze97a}), one can even give explicit
   values $d_k$ (resp.\ $d_k'$) for every $k>0$ such that
   $\pi^*L$ is $k$-jet ample (resp.\ $k-$very ample) whenever $d\ge d_k$
   (resp.\ $d\ge d'_k$).
\end{example}

   The previous example shows that the positivity can very well
   increase after taking a pullback.
   In general, however, it does not need to increase
   at all, as the following example shows:

\begin{example}\rm
   Consider a product of elliptic curves $X=E_1\times\dots\times E_n$
   and the product polarization
   $$
      L=\O_X\(\sum_{i=1}^n pr^*_i(0)\) \ ,
   $$
   where $pr_i$ is the
   $i$-th projection. The choice of a $d$-torsion point on $E_n$
   determines a covering
   $\tilde E_n\to E_n$ which in turn induces a cyclic covering
   $Y=E_1\times\dots\times E_{n-1}\times
   \tilde E_n\to X$ such that
   $$
      \pi^*L=\O_Y\(\sum_{i=1}^{n-1} pr^*_i(0)
      +d\cdot pr_n^*(0)\) \ .
   $$
   The line bundles $\O_X((k+2)L)$ and
   $\O_Y((k+2)\pi^*L)$ are then $k$-jet ample (by \cite{BauSze97}), but
   neither of them is $\ell$-jet ample (or even $\ell$-very ample) for
   any $\ell>k$, since the restricted bundles
   $\O_X((k+2)L)|E_n$ and
   $\O_Y((k+2)\pi^*L)|\tilde E_n$ are only of
   degree $k+2$ (cf.\ \cite[Proposition 2.1]{BelSom93a}).
\end{example}

  Even if $L$ is positive, the positivity of $\pi^*L$ may drop if the
  assumptions on the line bundles $L-qM$ are not satisfied. The following
  example illustrates their importance.

\begin{example}\rm
  Let $X\cong \bbP^n$ and let $r$ and $d$ be integers greater or equal $2$.
  Let $M=\calo_{\bbP^n}(r)$ and let $B$ be a smooth member of $|dM|$. Let
  $\pi:Y\lra X$ be the cyclic covering defined by $B$ and let
  $L=\calo_{\bbP^n}((d-1)r)$. Then $L$ is $d-$jet ample but $\pi^*L$
  fails to be that positive. The reason is that $L-qM$ fails to be
  $(d-q)-$jet ample for $q=d-1$. Indeed, as in the proof of the case $3$ in
  Theorem \ref{jet} let $y\in R$ be a point on the ramification divisor
  $R\subset Y$ and $u_1,\dots,u_n$ be local coordinates such that $R$
  is locally defined by $u_1=0$. From the decomposition
  $$H^0(\pi^*L)\cong H^0(\calo_{\bbP^n}((d-1)r))\oplus
    u_1\cdot H^0(\calo_{\bbP^n}((d-2)r))\oplus\dots\oplus
    u_1^{d-1}\cdot H^0(\calo_{\bbP^n})$$
  we infer that no section of $\pi^*L$ generates a jet of the form
  $u_1^{d-1}u_i$ for $i=2,\dots,n$.
  In fact Theorem \ref{jet} implies that $\pi^*L$ is $(d-1)-$ jet ample.
\end{example}
  
\begin{remark}\rm
  More generally, in the case of branched coverings,
  using the decomposition (\ref{decomposition}) one can
  show as in the previous example
  that if $L-qM$ fails to generate $(k-q)-$jets
  at the points of the branched locus, $\pi^*L$ fails
  to be $k-$jet ample.
\end{remark}

  Now, we want to show that the assumptions of
  Theorem \ref{very} cannot be weakened.

\begin{example}(Geiser involution)\ \rm  
  Let $B\subset X\cong\bbP^2$ be a smooth plane quartic and let
  $\pi:Y\lra X$ be the double covering branched over $B$.
  Then $Y$ is the Del Pezzo surface of degree $2$
  and $-K_Y\cong\pi^*\calo_{\bbP^2}(1)$, (see \cite[10.2.4]{BelSom95}).
  Theorem \ref{very} states that $-kK_Y$ is $k-$very ample
  for $k\geq 2$. This reproves a result of \cite{DiR96}.
  Moreover, in \cite{DiR96} it was shown that $-kK_Y$ is not
  $(k+1)-$very ample. This implies that our theorem is sharp.
\end{example}

  Finally, we show how our results can be applied in certain
  situations to study higher order embeddings of blown up
  varieties.

\begin{example}(Bertini involution)\ \rm
  Let $S$ be the Del Pezzo surface of degree $1$. Then $|-K_S|$
  has a base point $x$ and $|-2K_S|$ is globally generated and
  defines a $2:1$ mapping $p:S\lra Q$ onto a quadric cone $Q\subset\bbP^3$.
  It is known (cf. \cite[10.4.3]{BelSom95}) that $p$ is branched
  over a smooth curve and the vertex $v=p(x)$ of $Q$. Let
  $\tau:\bbF_2\lra Q$ be the blowing up of $Q$ at $v$ with the
  exceptional divisor $D$ and let
  $\sigma:X\lra S$ be the blowing up of $S$ at $x$ with the exceptional
  divisor $E$. Then the following diagram
  \begin{commdiag}
    X & \Rarrow^{\pi} & \bbF_2 \\
    \Darrow\llabel{\sigma} & & \Darrow\rlabel{\tau} \\
    S & \Rarrow^{p} & Q
  \end{commdiag}
  is commutative and $\pi$ is a double covering branched over a smooth
  divisor homologous to $2(2D+3f)$, where $f$ denotes the pull-back
  of the ruling on $Q$ by $\tau$. It is interesting to investigate the
  positivity of the line bundles of the form $-sK_X-tE$, for $s,t>0$. Notice that for $s,t\leq 0$ the line bundle is not even ample. We observe that
  $-sK_X-tE=s\sigma^*(-K_S)-(s+t)E$ but since the jet ampleness of $-sK_S$
  is not known for we cannot apply the results from \cite{BelSomJA}
  directly. On the other hand
  $-sK_X-tE=s\pi^*f-t\pi^{-1}(D)$.
  This suggests to look for the positivity of the line bundles
  coming from $\bbF_2$ via $\pi$. We recall 
  that a line bundle $L=aD+bf$ is $k-$jet ample on $\bbF_2$
  (equivalently $k-$very ample,
  cf. e.g. \cite{DiR2}, \cite{BelSom93a}) if and only if
  $L.D=-2a+b\geq k$ and $L.f=a\geq k$.
  Then Theorem \ref{jet} implies that $\pi^*L=-bK_X-2aE$ is $k-$jet ample
  if $a\geq k+1$ and $b\geq 3k$ and Theorem \ref{very} asserts
  that $\pi^*L$ is $k-$very ample for 
  $a\geq \rounddown{\frac{k-1}2}+2$ and
  $b\geq 3\rounddown{\frac{k-1}2}+3$.
\end{example}

%****************************************************************************

\small
\makeatletter
\def\@listi{\topsep=0cm\parsep=0cm\itemsep=0cm}

%****************************************************************************

\bigskip
\bigskip
\bigskip
\parindent=0cm

   Thomas Bauer,
   Mathematisches Institut,
   Universit\"at Erlangen-N\"urnberg,
   Bismarckstra{\ss}e $1\frac12$,
   D-91054 Erlangen,
   Germany

   E-mail: {\tt bauerth@mi.uni-erlangen.de}

\bigskip
   Sandra Di Rocco,
   KTH, Department of Mathematics,
   S-100 44 Stockholm,
   Sweden

   E-mail: {\tt sandra@math.kth.se}

   Current address: Max-Planck-Institut f\"ur Mathematik,
   Gottfried-Claren-Str.\ 26,
   D-53225 Bonn, Germany

\bigskip
   Tomasz Szemberg,
   Instytut Matematyczny PAN,
   c/o Instytut Matematyki,
   Uniwersytet Jagiello\'nski,
   Reymonta 4,
   PL-30-059 Krak\'ow,
   Poland

   E-mail: {\tt szemberg@im.uj.edu.pl}
%****************************************************************************

\end{document}